\documentclass{cocv}
\usepackage{graphicx}
\begin{document}
\title{A method for finding a solution to the nonsmooth differential inclusion of a special structure}\thanks{The main results of this paper (Sections \ref{sc6}, \ref{sc7}, \ref{sc8}, \ref{sc9}) were obtained in IPME RAS and supported by Russian Science Foundation (grant 20-71-10032).}
\author{Alexander Fominyh}\address{Institute of Problems in Mechanical Engineering, Russian Academy of Sciences,  61, Bolshoy pr. V.O., St. Petersburg, 199178, Russian Federation, \\  St.\,Petersburg University, 7--9, Universitetskaya nab., St.\,Petersburg, 199034, Russian Federation}
%
%
\begin{abstract} The paper explores the differential inclusion of a special form. It is supposed that the support function of the set in the right-hand side of an inclusion may contain the maximum of the finite number of continuously differentiable  (in phase coordinates) functions. It is required to find a trajectory that would satisfy the differential inclusion with the boundary conditions prescribed and simultaneously lie on the surface given. Such problems arise while practical modeling of discontinuous systems and in other applied problems. The initial problem is reduced to a variational one. It is proved that the resulting functional to be minimized is superdifferentiable. The necessary minimum conditions in terms of superdifferential are formulated. The superdifferential (or the steepest) descent method in a classical form is then applied in order to find stationary points of this functional. Herewith, the functional is constructed in a such a way that one can verify whether the stationary point constructed is indeed a global minimum point of the problem. The convergence of the method proposed is proved. The method constructed is illustrated by examples.  \end{abstract}
%
%
\maketitle
\section{Introduction}
Differential inclusions are a powerful instrument for modeling dynamical systems. In this paper the differential inclusions of some special structure are explored. More specific: the support function of the set in the right-hand side of a differential inclusion contains the maximum functions of the finite number of continuously differentiable (in phase coordinates) functions. Such differential inclusions arise from systems of differential equations with discontinuous right-hand sides (when considering such systems moving in a sliding mode) and from some other practical problems. So it is required to find a trajectory of such a differential inclusion which simultaneously would satisfy the boundary conditions and lie on the ``discontinuity'' surface. 

Note that the majority of methods in literature consider only differential inclusions with a free right endpoint and use some classical approaches as Euler and Runge-Kutta schemes, various finite differences methods etc. (see, e. g., \cite {Sandberg}, \cite{Bastien}, \cite{Beyn}, \cite{Lempio}, \cite{Veliov}). A survey of difference methods for differential inclusions can be found in \cite{Dontchev}. Note one paper \cite{Schilling} where an algorithm to solve boundary value problems for differential inclusions was constructed. Let us also give some references \cite {Cernea}, \cite{Pappas}, \cite{Zhu}, \cite{Aru} with optimality conditions in problems with differential inclusion. In these papers differential inclusions of a rather general form are considered. There are cases of phase constraints as well as nonsmooth and nonconvex ones. On the other hand, the works listed are more of theoretical significance and some results seem hard to be employed in practice.  

The apparatus used for solving the problem in the paper is based on the ideas developed in the author previous works. The initial problem is reduced to a variational one, and the resulting functional to be minimized occurs to be nonsmooth (but only superdifferentiable). So the methods of nondifferentiable optimization are required to explore this problem. In paper \cite{Fom_opt_lett} a method for solving the classical nonsmooth (but only subdifferential) variational problem is proposed based on the idea of considering phase trajectory and its derivative as independent variables (and taking the natural connection between these variables into account via the special penalty term). Reducing solving differential inclusions to a variational problem via corresponding support functions was carried out in papers \cite{Fominyh_1},  \cite{Fominyh_2}, \cite{Fominyh_3}. Some other nonsmooth problems of optimal control and variational calculus were considered via similar methods in papers \cite{FomDolg}, \cite{Fom_nonsmooth}. The paper also uses some ideas of V. F. Demyanov scientific school on nondifferentiable minimization. For example, in book \cite{demmal} a modified subdifferentiable (steepest) descent method for minimizing the maximum of the finite number of continuously differentiable functions in the finite dimensional space is justified. 

In the paper presented some basic facts from the control of discontinuous systems of differential equations and sliding modes theory is used (see Section \ref{sc2}). It is well known that when the system state is in a
sliding mode, new useful properties of the model are observed. For example such motions may be optimal
in the sense of some criterion in the optimal control theory. The sliding modes are also used in order to stabilize the system, as well as to get rid of unwanted external disturbances. In some cases the system motion in a sliding mode may be described via differential inclusions with the special structure, and this is a motivation of considering such differential inclusions in this paper.  One can get in touch with the basics required as well as with the statement problem motivation regarding sliding modes in works \cite{Filippov}, \cite{Utkin}, \cite{Levant}. 
\section{Basic definitions and notations} \label{sc2}
In this paper we use the following notations. Let $C_{n} [0, T]$ be the space of $n$-dimensional continuous on $[0, T]$ vector-functions. Let also $P_{n} [0, T]$ be the space of piecewise continuous and bounded on $[0, T]$ \linebreak $n$-dimensional vector-functions. We also require the space $L^2_n [0, T]$ of square-summable on $[0, T]$ $n$-dimensional vector-functions. If $X$ is some normed space, then $||\cdot||_X$ denotes its norm and $X^*$ denotes the space conjugate to the space given.

We will assume that each trajectory $x(t)$ is a piecewise continuously differentiable vector-function. Let $t_0 \in [0, T)$ be a point of nondifferentiability of the vector-function $x(t)$, then we assume that $\dot x(t_0)$ is a right-hand derivative of the vector-function $x(t)$ at the point~$t_0$ for definiteness. Similarly, we assume that $\dot x(T)$ is a left-hand derivative of the vector-function $x(t)$ at the point~$T$. Now with the assumptions and the notations taken we can suppose that the vector function $x(t)$ belongs to the space $C_{n} [0, T]$ and that the vector function~$\dot x(t)$ belongs to the space $P_{n} [0, T]$.

For some arbitrary set $F \subset R^n$ define the support function of the vector $\psi~\in~R^n$ as $c(F, \psi) = \sup \limits_{f\in F}\langle f, \psi \rangle$ where $\langle a, b \rangle$ is a scalar product of the vectors $a, b \in R^n$. Denote $S_n$ a unit sphere in $R^n$ with the center in the origin, also let $B_r(c)$ be a ball with the radius $r \in R$ and the center $c \in R^n$. Let the vectors $\bf{e_i}$, $i =\overline{1,n}$, form the standard basis in $R^n$. The sum $P+Q$ of the sets $P, Q \subset R^n$
is their Minkowski sum, while $\lambda P$ with $\lambda \in R$ is the Minkowski product; if $q$ is some vector from $R^n$, then $Pq = \{ \langle p, q \rangle \ | \ p \in P\}$. Let $0_n$ denote a zero element of a functional space of some $n$-dimensional vector-functions and $\bf{0_n}$ denote a zero element of the space $R^n$. Let $E_m$ be an identity matrix and ${\bf {O}_{n, m}}$ --- a zero matrix in the space $R^n \times R^m$. If $\varphi(x) = \min\limits_{i = \overline{1, M}} f_i(x)$, where $f_i(x): R^n \rightarrow R$, $i = \overline{1, M}$, are some functions, then we call the function $f_{\overline i}(x)$, $\overline i \in \{1,\dots,M\}$, an active one at the point $x_0 \in R^n$, if $\overline i \in R(x_0) = \left\{i \in \{1,\dots,M\right\} \ | \ f_i(x_0) = \varphi(x_0) \} $.


In the paper we will use both superdifferentials of functions in a finite-dimensional space and superdifferentials of functionals in a functional space. Despite the fact that the second concept generalizes the first one, for convenience we separately introduce definitions for both of these cases and for those specific functions (functionals) and their variables and spaces which are considered in the paper. 

Consider the space $R^n \times R^n$ with the standard norm. Let $d = [d_1, d_2] \in R^n \times R^n$ be an arbitrary vector. Suppose that at the point $(x, z)$ there exists such a convex compact set $\overline \partial h(x,z)$ $\subset R^n \times R^n$ that 
\begin{equation}
\label{0.1}
\frac{\partial h(x,z)}{\partial d} = \lim_{\alpha \downarrow 0} \frac{1}{\alpha} \big(h(x+\alpha d_1, z + \alpha d_2) - h(x,z)\big) = \min_{w \in \overline \partial h(x,z)} \langle w, d \rangle. 
\end{equation}

In this case the function $h(x,z)$ is called superdifferentiable at the point $(x, z)$, and the set $\overline \partial h(x,z)$ is called the superdifferential of the function $h(x,z)$ at the point $(x,z)$.

From expression (\ref{0.1}) one can see that the following formula
$$ h(x + \alpha d_1, z + \alpha d_2) = h(x, z) + \alpha  \frac{\partial h(x,z)}{\partial d} + o(\alpha, x, z, d), $$
$$\quad  \frac{o(\alpha, x, z, d)}{\alpha} \rightarrow 0, \ \alpha \downarrow 0,$$
holds true. 

%

If the function $\varsigma(\xi)$ is differentiable at the point $\xi_0 \in R^\ell$, then its superdifferential at this point is represented in the form \begin{equation}
\label{0.4}\overline \partial \varsigma(\xi_0) = \{ \varsigma'(\xi_0) \} 
\end{equation}
 where $\varsigma'(\xi_0)$ is a gradient of the function $\varsigma(\xi)$ at the point $\xi_0$.
Note also that the superdifferential of the finite sum of superdifferentiable functions is the sum of the superdifferentials of summands, i. e. if the functions $\varsigma_k(\xi)$, $k = \overline{1, r}$, are superdifferentiable at the point $\xi_0 \in R^\ell$, then the function $\varsigma(\xi) = \sum_{k=1}^r \varsigma_k(\xi)$ superdifferential at this point is calculated by the formula 
\begin{equation}
\label{0.5} \overline \partial \varsigma(\xi_0) = \sum_{k=1}^r \overline\partial \varsigma_k(\xi_0).
\end{equation}

Consider the space $C_n[0, T] \times P_n[0, T]$ with the norm $L_n^2 [0, T] \times L_n^2 [0, T]$. Let $g = [g_1, g_2] \in C_n[0, T] \times P_n[0, T]$ be an arbitrary vector-function. Suppose that at the point $(x, z)$ there exists such a convex weakly$^*$ compact set $\overline \partial {I(x, z)} \subset \big( C_n[0, T] \times P_n[0, T],$ $|| \cdot ||_{L_n^2 [0, T] \times L_n^2 [0, T]} \big) ^*$ that  
\begin{equation}
\label{0.6} 
\frac{\partial I(x, z)}{\partial g} = \lim_{\alpha \downarrow 0} \frac{1}{\alpha} \big(I(x+\alpha g_1, z+\alpha g_2) - I(x, z)\big) = \min_{w \in \overline \partial I(x, z)} w(g).
\end{equation}

In this case the functional $I(x,z)$ is called superdifferentiable at the point $(x, z)$ and the set $\overline \partial {I(x, z)} $ 
is called a superdifferential of the functional $I(x, z)$ at the point $(x,z)$. 

From expression (\ref{0.6}) one can see that the following formula 
$$ 
I(x + \alpha g_1, z + \alpha g_2) = I(x, z) + \alpha  \frac{\partial I(x,z)}{\partial g} + o(\alpha, x, z, g),
$$
$$ \quad  \frac{o(\alpha, x, z, g)}{\alpha} \rightarrow 0, \ \alpha \downarrow 0,$$
holds true.

%
%

\section{Statement of the problem}
Consider the differential inclusion
\begin{equation}
\label{1}
\dot x_i \in A_i x + [\underline a_i, \overline a_i] |x| [-1, 1] = A_i x + [-\overline a_i, \overline a_i] |x|, \quad i = \overline {1,n},
\end{equation}
with the initial point
\begin{equation}
\label{2}
x(0) = x_{0}
\end{equation}
and with the desired endpoint
\begin{equation}
\label{3}
x_j(T) = {x_T}_j, \quad j \in J,
\end{equation}
In formula (\ref{1}) $A_i$ is the $i$-th row of the constant $n \times n$ matrix $A$, $i = \overline{1,n}$, and $\underline a_i$, $\overline a_i$, $i = \overline{1,n}$, are given nonnegative numbers and $\underline a_i = \overline a_i =0$, $i =\overline{m+1, n}$. The system is considered on the given finite time interval $[0, T]$. We assume that $x(t)$ is an $n$-dimensional continuous vector-function of phase coordinates with a piecewise-continuous and bounded derivative on the segment $[0, T]$. In formula (\ref{2}) $x_0 \in R^n$ is a given vector; in formula (\ref{3}) ${x_T}_j$ are given numbers corresponding to those coordinates of the state vector which are fixed at the right endpoint, here $J \subset \{1, \dots, n\}$ is a given index set.  

Let $F(x) = (f_1(x), \dots, f_n(x))'$ where $f_1(x), \dots, f_n(x)$ run through the corresponding sets $F_1(x), \dots F_n(x)$ from the right-hand sides of inclusions (\ref{1}), so we can rewrite the given inclusions in the form 
\begin{equation}
\label{4}
\dot{x} \in F(x).
\end{equation}

Let us also take into consideration the surface 
\begin{equation}
\label{5}
s(x) = {\bf 0_m}
\end{equation}
where $s(x)$ is a known continuously differentiable $m$-dimensional vector-function. 


 We formulate the problem as follows: it is required to find such a trajectory $x^{*} \in C_{n}[0, T]$ (with the derivative $\dot x^{*} \in P_{n}[0, T]$) which moves along surface (\ref{5}), satisfies differential inclusion (\ref{4}) while $t \in [0,T]$ and meets boundary conditions (\ref{2}), (\ref{3}). Assume that there exists such a solution.

Let us discuss one practical problem which leads to the statement of the problem above. Consider the system
\begin{equation}
\label{6}
\dot{x} = A x + B u
\end{equation}
with boundary conditions (\ref{2}), (\ref{3}).
In formula (\ref{6}) $A$ is a constant $n \times n$ matrix, $B$ is a constant $n \times m$ matrix. For simplicity we suppose that $B = [E_m, {\bf{O}_{n-m, m}}]$. The system is considered on the given finite time interval $[-t^*, T]$ (here $T$ is a given final time moment; see comments on the time moment $t^*$ below). We suppose that $x(t)$ is an $n$-dimensional continuous vector-function of phase coordinates with a piecewise-continuous and bounded derivative on $[-t^*, T]$; the structure of the $m$-dimensional control $u$ will be specified below.

Let also ``discontinuity'' surface (\ref{5}) be given.

Consider the following form of controls:
\begin{equation}
\label{7}
u_i = -\alpha_i |x| \mathrm{sign}(s_i(x)), \quad i = \overline{1, m}
\end{equation} 
where $\alpha_i \in [\underline {a}_i, \overline {a}_i]$, $i = \overline{1, m}$, are some positive numbers which are sometimes called gain factors.

In book \cite{Utkin} it is shown that if surface (\ref{5}) is a hyperplane, then under natural assumptions and with sufficiently big values of the factors $\alpha_i$, $i = \overline{1,m}$, controls (\ref{7}) ensure system (\ref{6}) hitting a small vicinity of this surface $s(x) = {\bf 0_m}$ from arbitrary initial state (\ref{2}) in the finite time $t^*$ and further staying in this neighborhood with the fulfillment of the condition   $s_i(x(t)) \rightarrow 0$, $i = \overline{1,m}$, at $t \rightarrow \infty$, i. e. controls (\ref{7}) ensure the stability ``in big'' of the system (\ref{6}) sliding mode. In \cite{Utkin} one may also find the estimates on the time moment $t^*$. Here we assume that all the conditions required are already met, and the numbers $\underline a_i$, $\overline a_i$, $i = \overline{1, m}$, are taken sufficiently big. Now we are interested in the behavior of the system on the ``discontinuity'' surface (on the time interval $[0,T]$). 
 


We see that the right-hand sides of the first $m$ differential equations in system (\ref{6}) with controls  (\ref{7}) are discontinuous on the surfaces $s_i(x) = 0$, $i = \overline{1,m}$. So one has to use one of the known definitions of a discontinuous system solution. 

Let us use one of the classical variants of such a definition \cite{Filippov}, \cite{AizerPyatnizkyi}. The essence of this definition is as follows.
 On the finite time interval $[0, T]$ consider the system $\dot x = f(x, u_1(x, t), \dots, u_m(x, t), t)$, 
 in which the vector-function $f(x, u_1, \dots, u_m, t)$ is continuous in all its arguments, and the vector-functions $u_i (x, t)$, $i = \overline{1,m}$, 
 are discontinuous on the sets $s_i(x) = 0$, $i = \overline{1,m}$, respectively. 
 At every point $(x, t)$ of discontinuity of the vector-function $u_i (x, t)$, $i = \overline{1,m}$, 
 a closed set $U_i (x, t)$, $i = \overline{1,m}$, must be defined. It is a set of possible values of the variable  
 $u_i$ of the function $f(x, u_1, \dots, u_m, t)$. 
 Denote $F (x, t) = f(x, u_1, ..., u_m, t)$  
the set of the function $f(x, u_1, ..., u_m, t)$ values at the fixed variables 
 $x, t$, and while $u_1, \dots, u_m$ run through the sets 
 $U_1(x, t), \dots, U_m(x, t)$ respectively. Then the solutions of this differential inclusion are taken as solutions of the original differential equation with a discontinuous right-hand side. 
In physical systems the sets $U_i (x, t)$, $i = \overline{1,m}$, usually correspond to different blocks and are assumed to be convex. At every point $(x, t)$ of discontinuity of the vector-function $u_i (x, t)$, $i = \overline{1,m}$, the set $U_i (x, t)$ must also contain all the limit points of all sequences $v_k \in U_i (x_k, t_k)$, where $x_k \rightarrow x$ and $t_k \rightarrow t$ if $k \rightarrow \infty$. If the control is of form (\ref{7}), then it is natural to consider $U_i(x, t) = \mathrm{co} \{-\overline a_i |x|, \overline a_i |x|\}$, $i = \overline{1,m}$, as such sets. 

Thus, accordingly to this definition a solution of the discontinuous system considered satisfies differential inclusion (\ref{1}) above. Note that the more detailed version of the definition given is presented in \cite{AizerPyatnizkyi}. It has a strict and rather complicated form so we don't consider the details here. For our purposes it is sufficient to postulate that the system (\ref{6}) solutions with controls (\ref{7}) employed are the solutions of inclusion (\ref{1}) by definition. 

\begin{rmrk}
Instead of trajectories from the space $C_n [0, T]$ with derivatives from the space $P_n [0, T]$ one may consider absolutely continuous trajectories on the interval $[0, T]$ with measurable and almost everywhere bounded derivatives on $[0, T]$ what is more natural for differential inclusions. The choice of the solution space in the paper is explained by the possibility of its practical construction.
\end{rmrk}

\section{Reduction to a variational problem} We will sometimes write $F$ instead of $F(x)$ for brevity. 
Insofar as  $\forall x \in R^n$ the set $F(x)$ 
is a convex compact set in $R^n$, then inclusion (\ref{4}) 
may be rewritten as follows {\cite{Blagodatskih}}:
$$
\dot x_i(t) \psi_i(t) \leq c(F_i(x(t)), \psi_i(t)) \quad \forall \psi_i(t) \in S_1, \quad \forall t \in [0, T], \quad i = \overline{1,n}.
$$

Calculate the support function of the set $F_i$.  For this note that the set $F_i$ is a one-dimensional ``ball'' with the the center 
$$ c_i(x) = A_i x , \quad i = \overline{1,n},$$
and with the ``radius'' 
$$ r_i(x) = \overline a_i |x|, \quad i = \overline {1,m},$$
$$ r_i(x) = 0, \quad i = \overline{m+1, n}.$$

So the support function of the set $F_i$ can be expressed \cite{Blagodatskih} by the formula
$$ c(F_i(x),\psi_i) = \psi_i A_i x + \overline a_i |x| |\psi_i|, \quad i = \overline {1,m},$$
$$ c(F_i(x),\psi_i) = \psi_i A_i x, \quad i = \overline{m+1, n}.$$

%
 
We see that the support function of the set $F_i$ is continuously differentiable in the phase coordinates $x$ if $x_i \neq 0$, $i = \overline{1, n}$. 

Denote $z(t) = \dot x(t)$, $z \in P_n [0, T]$, then from (\ref{2}) one has
\begin{equation}
\label{3.11''}
x(t) = x_0 + \int_0^t z(\tau) d\tau.
\end{equation}

Let us now realize the following idea. ``Forcibly'' consider the points $z$ and $x$ to be ``independent'' variables. Since, in fact, there is relationship (\ref{3.11''}) above between these variables (which naturally means that the vector-function $z(t)$ is a derivative of the vector-function $x(t)$), let us take this restriction into account by using the functional
 \begin{equation} 
\label{3.17}
 \upsilon(x, z) = \frac{1}{2} \int_0^T \left( x(t) - x_0 - \int_0^t z(\tau) d \tau \right)^2 dt.
\end{equation} 
Besides it is seen that condition (\ref{2}) on the left endpoint is also satisfied if $v(x, z) = 0$. 

For $i = \overline{1,n}$ put
$$
\ell_i(\psi_i, x, z) = \langle z_i, \psi_i \rangle - c ( F_i(x), \psi_i ),
$$
$$
h_i(x, z) = \max_{\psi_i \in S_1} \max \{ 0, \ell_i(\psi_i, x, z) \},
$$
then put
$$h(x,z) = (h_1(x,z), \dots, h_n(x,z))'$$
and construct the functional
\begin{equation}
\label{3.13}
\varphi(x, z) = \frac{1}{2} \int_0^T h^2 \big( x(t), z(t) \big) dt.
\end{equation}

It is not difficult to check that for functional (\ref{3.13}) the relation
$$
\left\{
\begin{array}{ll}
\varphi(x, z) = 0, \ &\text{if} \ \dot x_i(t) \psi_i(t)  \leq c(F_i(x(t)), \psi_i(t)) \quad \forall \psi_i(t) \in S_1, \quad \forall t \in [0, T], \quad i = \overline{1,n}. \\
\varphi(x, z) > 0, \ &\text{otherwise},
\end{array}
\right.
$$
holds true, i. e. inclusion (\ref{4}) takes place iff $\varphi(x, z) = 0$.

Introduce the functional
\begin{equation}
\label{3.14}
\chi(z) = \frac{1}{2} \sum_{j \in J} \left( {x_0}_j + \int_0^T z_j(t) dt - {x_T}_j \right)^2.
\end{equation}

 We see that if $v(x, z) = 0$, then condition (\ref{3}) on the right endpoint is satisfied iff $\chi(z) = 0$.

Introduce the functional 
\begin{equation}
\label{3.14'}
\omega(x) = \frac{1}{2} \int_0^T s^2 \big( x(t) \big ) dt.
\end{equation}

Obviously, the trajectory $x(t)$ belongs to surface (\ref{5}) at each $t \in [0, T]$ iff $\omega(x) = 0$.

Finally construct the functional
\begin{equation} 
\label{3.16} 
I(x, z) = \varphi(x, z) + \chi(z) + \omega(x) + \upsilon(x, z).
\end{equation}

So the original problem has been reduced to minimizing functional (\ref{3.16}) on the space $C_n[0, T] \times P_n [0, T]$. Denote $z^*$ a global minimizer of this functional. Then $$ x^*(t) = x_0 + \int_0^t z^*(\tau) d\tau $$
is a solution of the initial problem. 

\begin{rmrk}
The structure of the functional $\varphi(x, z)$ is natural as the value $ h_i(x(t), z(t)) $, $i = \overline{1, n}$,
at each fixed $t \in [0, T]$ is the Euclidean distance from the point $z_i(t)$
to the set $F_i (x(t))$; functional (\ref{3.13}) is half the sum of squares of the deviations in $L^2_n[0,T]$ norm of the trajectories $z_i(t)$ from the sets $F_i(x)$, $i = \overline{1, n}$, respectively. The meaning of functionals (\ref{3.17}), (\ref{3.14}), (\ref{3.14'}) structures is obvious.
\end{rmrk}

Despite the fact that the dimension of functional $I(x, z)$ arguments is $n$ more the dimension of the initial problem (i. e. the dimension of the desired point $x^*$), the structure of its superdifferential (in the space $C_n[0, T] \times P_n[0, T]$ as a normed space with the norm $L_n^2 [0, T] \times L_n^2 [0, T]$), as will be seen from what follows, has a rather simple form. This fact will allow us to construct
a numerical method for solving the original problem. 

\section{Minimum conditions of the functional $I(x, z)$}
In this section referring to superdifferential calculus rules (\ref{0.4}), (\ref{0.5}), we mean their known analogues in a functional space \cite{Dolgopolikcodiff}. 

Using classical variation it is easy to prove the Gateaux differentiability of the functional $\chi(z)$, we have
$$
 \nabla \chi(z) = \sum_{j \in J} \left( {x_0}_j + \int_0^T z_j(t) dt - {x_T}_j \right) {\bf e_j}.  
 $$
 By superdifferential calculus rule (\ref{0.4}) one may put 
\begin{equation} \label{4.44} \overline \partial \chi(z) =  \left\{ \sum_{j \in J} \left( {x_0}_j + \int_0^T z_j(t) dt - {x_T}_j \right) {\bf e_j} \right\}. \end{equation}
 
Using classical variation it is easy to prove the Gateaux differentiability of the functional $\omega(x)$, we have
$$
 \nabla \omega(x) = \sum_{i=1}^m s_i\big(x(t)\big) \frac{\partial s_i\big(x(t)\big)} {\partial x}.  
 $$
 By superdifferential calculus rule (\ref{0.4}) one may put 
\begin{equation} \label{4.6} \overline \partial \omega(x)  =  \left\{ \sum_{i=1}^m s_i\big(x(t)\big) \frac{\partial s_i\big(x(t)\big)} {\partial x} \right\}. \end{equation}

 Using classical variation and integration by parts it is also not difficult to check the Gateaux differentiability of the functional $\upsilon(x, z)$, we obtain
 $$\nabla \upsilon(x, z, t) = \begin{pmatrix}
\displaystyle{ x(t) - x_0 - \int_0^t z(\tau) d\tau } \\
\displaystyle{ -\int_t^T \left( x(\tau) - x_0 - \int_0^\tau z(s) ds \right) d\tau  }
\end{pmatrix}. 
$$   
 
  By superdifferential calculus rule (\ref{0.4}) one may put 
\begin{equation} \label{4.55} \overline \partial \upsilon(x, z)  =  \left\{ \begin{pmatrix}
\displaystyle{ x(t) - x_0 - \int_0^t z(\tau) d\tau } \\
\displaystyle{ -\int_t^T \left( x(\tau) - x_0 - \int_0^\tau z(s) ds \right) d\tau  }
\end{pmatrix} \right\}. \end{equation}
 
 Explore the differential properties of the functional $\varphi(x,z)$. For this, first give the following formulas for calculating the superdifferential $\overline \partial h^2(x, z)$ at the point $(x, z)$. At $i = \overline{1,m}$ one has
\begin{equation}
\label{4.5}
\overline \partial \left( {\textstyle \frac{1}{2}} \, h_i^2(x, z) \right) = h_i(x,z)\bigg(\psi^*_i {\bf e_{i+n}} - \psi^*_i [A'_i, {\bf 0_n}] + \sum_{j=1}^n \overline \partial \left(-\overline{a}_i |x_j| |\psi_i^*| \right)\bigg),
\end{equation}
where at $j = \overline{1,n}$ we have
$$ \overline \partial(-\overline{a}_i |x_j| |\psi_i^*|) = \left\{
\begin{array}{lll}
-\overline{a}_i |\psi_i^*| {\bf e_j}, \ &\text{if} \ x_j > 0,  \\
\overline{a}_i |\psi_i^*| {\bf e_j}, \ &\text{if} \ x_j < 0, \\
\mathrm{co} \big \{-\overline{a}_i |\psi_i^*| {\bf e_j},  \overline{a}_i |\psi_i^*| {\bf e_j} \big \}, \ &\text{if} \ x_j = 0.
\end{array}
\right.
$$
At $i = \overline{m+1,n}$ one has
$$\overline \partial \left( {\textstyle \frac{1}{2}} \, h_i^2(x, z) \right) = h_i(x,z)\big(\psi^*_i {\bf e_{i+n}} - \psi^*_i [A'_i, {\bf 0_n}]\big).$$ In the formulas given the value $\psi^*_i$ is such that $\max \{ 0, \ell_i(\psi_i^*(x, z), x, z) \} = \max_{\psi_i \in S_1} \max \{ 0, \ell_i(\psi_i, x, z) \}$, \linebreak $i = \overline{1, n}$. In \cite{FomDolg} it is shown that if $h_i(x, z) > 0$ then $\psi_i^*(x,z)$ is unique and continuous in $(x,z)$. In \cite{Fominyh_sl_m} it is shown that the functional $\varphi(x,z)$ is superdifferentiable and that its superdifferential is determined by the corresponding integrand superdifferential. 
 
\begin{thrm}
Let the interval $[0, T]$ may be divided into a finite number of intervals, in every of which each phase trajectory is either identically equal to zero or retains a certain sign. Then the functional
 $
 \varphi(x, z) 
 $
is superdifferentiable, i. e. 
$$
\frac{\partial \varphi(x, z)}{\partial g} = \lim_{\alpha \downarrow 0} \frac{1}{\alpha} \big(\varphi(x+\alpha g_1, z+\alpha g_2) - \varphi(x, z)\big) = \min_{w \in \overline \partial \varphi(x, z)} \int_0^T \langle w(t), g(t) \rangle dt, 
$$
where $g = [g_1, g_2] \in C_n[0,T] \times P_n[0,T]$ and the set $\overline \partial \varphi(x, z)$ is defined as follows
\begin{equation}
\label{4.2} 
\overline \partial \varphi(x, z) = \Big\{ w = [w_1, w_2] \in L_{n}^\infty[0, T] \times L_{n}^\infty[0, T] \ \big|
 \end{equation}
$$ [w_1(t), w_2(t)] \in \overline \partial \left( {\textstyle \frac{1}{2}} \, h^2(x(t), z(t)) \right) \quad for \ a. e. \ t \in [0, T] \Big\}.$$
\end{thrm}

Using formulas (\ref{3.16}), (\ref{4.44}), (\ref{4.6}), (\ref{4.55}), (\ref{4.2}) obtained and superdifferential calculus rule (\ref{0.5}) we have the following final formula for calculating the superdifferential of the functional $I(x, z)$ at the point $(x, z)$
\begin{equation}
\label{4.3}
\overline \partial I(x, z) = \overline \partial \varphi(x,z) + \overline \partial \overline \chi(x,z) + \overline \partial \overline \omega(x,z) + \overline \partial \upsilon(x,z),
\end{equation} 
where formally $\overline \chi(x, z) := \chi(z)$, $\overline \omega(x, z) := \omega(x)$.

Using the known minimum condition \cite{Dolgopolikcodiff} of the functional $I(x, z)$ at the point $(x^*, z^{*})$ in terms of superdifferential, we conclude that the following theorem is true.

\begin{thrm}\label{th2}
Let the interval $[0, T]$ may be divided into a finite number of intervals, in every of which each phase trajectory $x^*_i(t)$, $i=\overline{1,n}$, is either identically equal to zero or retains a certain sign. In order for the point $(x^{*}, z^*)$ to minimize the functional $I(x, z)$, it is necessary to have 
\begin{equation}
\label{4.4}
\overline \partial I (x^*(t), z^{*}(t)) = \{\bf 0_{2n}\}
\end{equation}
at each $t \in [0, T]$ where the expression for the superdifferential $\overline \partial I (x, z)$ is given by formula (\ref{4.3}). If one has $I(x^*, z^{*})$, then condition (\ref{4.4}) is also sufficient.
\end{thrm}

Theorem \ref{th2} already contains a constructive minimum condition, since on its basis it is possible to construct the steepest (the superdifferential) descent direction, and for solving each of the subproblems arising during this construction there are known efficient algorithms for solving them (see Section \ref{sc6} below). Once the steepest descent direction is constructed, one is able to use it in order to apply some of nonsmooth optimization methods (see Section \ref{sc7} below). Note once again that without ``separation'' of the variables $x$ and $z$ proposed the considered functional superdifferential would have a very complicated structure. That would make constructing the superdifferential descent direction of this functional a significantly difficult (and as it seems, practically impossible) problem. 

\section{A more general case of the set $F(x)$ support function}\label{sc6}
Consider now a more general case when the support functions of the corresponding sets $F_i(x)$, $i = \overline{1,n}$, are of the form 
$$ c(F_i(x),\psi_i) = \psi_i A_i x + \sum_{j=1}^{r} \overline a_{ij} \max\big\{f_{i, j_1}(x) \psi_i, \dots, f_{i, j_{k(j)}}(x) \psi_i\big\},$$
where $f_{i, j_1}(x), \dots,  f_{i, j_{k(j)}}(x)$, $i = \overline{1, n}$, $j = \overline{1, r}$ (for simplicity of presentation we
suppose that $r$ is the same for each $i = \overline{1, n}$), are continuously differentiable functions and $\overline a_{ij}$, $i = \overline{1, n}$, $j = \overline{1, r}$, are some nonnegative numbers. (We also still consider only convex and compact sets $F_1(x), \dots F_n(x)$ at every $x \in R^n$).

This case indeed is more general as we have $|x_i| |\psi_i| = |x_i \psi_i| = \max\{x_i \psi_i, -x_i \psi_i\}$.   

One can also give a practical problem which leads to such a system. Let from some physical considerations the ``velocity'' $\dot x_1$ of an object lie in the range $[\min\{x_1, x_2, x_3\}, \max\{x_1, x_2, x_3\}]$ of the ``coordinates'' $x_1$, $x_2$, $x_3$. The segment given may be written down as ${\mathrm{co} \{x_1, x_2, x_3\} = \mathrm{co}\bigcup\limits_{i=1}^3 \{x_i\}}$. The support function of this set is \cite{Blagodatskih} $\max  \{x_1 \psi_1, x_2 \psi_1, x_3 \psi_1\}$.  

For simplicity consider the case $n=2$, $r = 1$, $k(1)=2$, and only the functions $\ell_1(\psi_1, x, z)$ and $h_1(x, z)$  (here we denote them $\ell(\psi_1, x_1, x_2, z_1)$ and $h(x_1, x_2, z_1)$ respectively) and the time interval $[t_1, t_2] \subset [0, T]$ of nonzero length; the general case is considered in a similar way. Then we have $ \ell (\psi_1, x_1, x_2, z_1) = z_1 \psi_1 -  a_1 x_1 \psi_1 -$ $- a_2 x_2 \psi_1 - b \max\{f_1(x) \psi_1, f_2(x) \psi_1\}$, where $a_1 := a_{1,1}$, $a_2 := a_{1,2}$, $b := \overline a_{1,1} > 0$, $f_1 := f_{1,1_1}$, $f_2 := f_{1,1_2}$. Fix some point $(x_1, x_2) \in R^2$. Let $f_1(x(t)) \psi_1(t) = f_2(x(t)) \psi_1(t)$ at $t \in [t_1, t_2]$; other cases may be studied in a completely analogous fashion.

a) Suppose that $h_1 (x, z) > 0$, i. e. $h_1(x, z) = \max_{\psi_1 \in S_1}\ell_1 (\psi_1, x, z) > 0$. 

Our aim is to apply the corresponding theorem on a directional differentiability from \cite{BonnansShapiro}. The theorem of this book considers the inf-functions so we will apply this theorem to the function $-\ell(\psi_1, x_1, x_2, z_1)$. For this, check that the function $h(x_1, x_2, z_1)$ satisfies the following conditions: \newline
i) the function $\ell (\psi_1, x_1, x_2, z_1)$ is continuous on $S_1 \times R^2 \times R$; \newline
ii)  there exist a number $\beta$ and a compact set $C \in R$ such that for every $(\overline x_1, \overline x_2, \overline z_1)'$ in the vicinity of the point $(x_1, x_2, z_1)'$ the level set 
$${lev}_\beta (-\ell(\cdot, \overline x_1, \overline x_2, \overline z_1)) = \{\psi_1 \in S_1 \ | \ -\ell(\psi_1, \overline x_1, \overline x_2, \overline z_1) \leq \beta\} $$
is nonempty and is contained in the set $C$; \newline
iii) for any fixed $\psi_1 \in S_1$ the function $\ell(\psi_1, \cdot, \cdot,\cdot)$ is directionally differentiable at the point $(x_1, x_2, z_1)'$; \newline
iv) if $d = [d_1, d_2] \in R^2 \times R$, $\gamma_n \downarrow 0$ and $\psi_{1_{n}}$ is a sequence in $C$, then $\psi_{1_{n}}$ has a limit point $\overline \psi_1$ such that 
$$\mathrm{lim} \sup\limits_{n \rightarrow \infty} \frac{-\ell(\psi_{1_{n}}, x + \gamma_n d_{1},  z+\gamma_n d_2) - (-\ell(\psi_{1_{n}}, x, z))}{\gamma_n} \geq \frac{\partial(-\ell(\overline \psi_1, x, z))}{\partial d}$$
where $\displaystyle{\frac{\partial\ell(\overline \psi_1, x_1, x_2, z_1)}{\partial d}}$ is the derivative of the function $\ell (\overline \psi_1, x_1, x_2, z_1)$ at the point $(x_1,x_2,z_1)'$ in direction~$d$. 

The verification of conditions i), ii) is obvious. 

In order to verify the condition iii), it is sufficient to observe that since $b > 0$, then for the fixed $\psi_1 \in S_1$ the function $- b \max\{f_1(x) \psi_1, f_2(x) \psi_1\}$ is superdifferentiable \cite{demvas} (hence, it is differentiable in directions) at the point $(x_1, x_2)'$; herewith, its superdifferential at this point is $ b \ \mathrm{co} \left \{ \left(
\psi_1 \frac{\partial f_1(x)}{\partial x_1},
\psi_1 \frac{\partial f_1(x)}{\partial x_2} \right)', \ 
\left(
\psi_1 \frac{\partial f_2(x)}{\partial x_1},
\psi_1 \frac{\partial f_2(x)}{\partial x_2} \right)'
\right \}$. 
An explicit expression of this function derivative at the point $(x_1, x_2)'$ in the direction $d_{1}$ is \linebreak $-b \max\left\{ \left \langle \psi_1 \frac{\partial f_1(x)} {\partial x}, d_1 \right \rangle, \left\langle \psi_1 \frac{\partial f_2(x)} {\partial x}, d_1 \right \rangle \right\} $. 

Finally, check condition iv). Let $[d_1, d_2] \in R^2 \times R$, $\gamma_n \downarrow 0$ and $\psi_{1_n}$ is some sequence from the set $C$. Calculate
$$\mathrm{lim} \sup\limits_{n \rightarrow \infty} \frac{-\ell(\psi_{1_n}, x_1 + \gamma_n d_{1,1}, x_2 + \gamma_n d_{1,2}, z_1 + \gamma_n d_2) - (-\ell(\psi_{1_n}, x_1, x_2, z_1))}{\gamma_n} =$$  $$=\mathrm{lim} \sup\limits_{n \rightarrow \infty} \frac{1}{\gamma_n} \Big(-(z_1 + \gamma_n d_2) \psi_{1_n} + a_1(x_1+\gamma_n d_{1,1}) \psi_{1_n}  + a_2(x_2+\gamma_n d_{1,2}) \psi_{1_n} +$$ $$+ b \max\big\{f_1(x_1+\gamma_n d_{1,1}, x_2 + \gamma_n d_{1,2}) \psi_{1_n}, f_2(x_1+\gamma_n d_{1,1}, x_2 + \gamma_n d_{1,2}) \psi_{1_n}\big\} + $$ $$+z_1 \psi_{1_n} - a_1 x_1 \psi_{1_n} - a_2 x_2 \psi_{1_n} - b \max\big\{f_1(x_1, x_2) \psi_{1_n}, f_2(x_1, x_2) \psi_{1_n}\big\} \Big)   = $$

$$= \mathrm{lim} \sup\limits_{n \rightarrow \infty} \frac{1}{\gamma_n} \bigg(-\gamma_n d_2 \psi_{1_n} + \gamma_n a_1 d_{1,1} \psi_{1_n} + \gamma_n a_2 d_{1,2} \psi_{1_n} + $$ $$+ b \max \left \{ \left [ f_1(x) + \gamma_n \bigg\langle \frac{\partial f_1(x)} {\partial x}, d_1 \bigg \rangle + o_1(\gamma_n, x, d)   \right] \psi_{1_n}, \left [ f_2(x) + \gamma_n \left\langle \frac{\partial f_2(x)} {\partial x}, d_1 \right \rangle + o_2(\gamma_n, x, d)   \right] \psi_{1_n} \right\} -$$ $$- b \max\big\{f_1(x_1, x_2) \psi_{1_n}, f_2(x_1, x_2) \psi_{1_n}\big\}  \bigg) \geq$$
$$ \geq \mathrm{lim} \sup\limits_{n \rightarrow \infty} \frac{1}{\gamma_n} \bigg(-\gamma_n d_2 \psi_{1_n} + \gamma_n a_1 d_{1,1} \psi_{1_n} + \gamma_n a_2 d_{1,2} \psi_{1_n} + $$ $$+ \gamma_n b \max \left \{  \bigg\langle \frac{\partial f_1(x)} {\partial x}, d_1 \bigg \rangle  \psi_{1_n}, \left\langle \frac{\partial f_2(x)} {\partial x}, d_1 \right \rangle \psi_{1_n} \right\} + b \min\big\{o_1(\gamma_n, x, d) \psi_{1_n}, o_2(\gamma_n, x, d) \psi_{1_n}\big\} \bigg)$$
where ${o_q(\gamma_n, x, d)}/{\gamma_n} \rightarrow 0, \ \gamma_n \downarrow 0$, $q=\overline{1,2}$, and the last inequality follows from the assumption \linebreak $f_1(x) \psi_1 = f_2(x) \psi_1$ (on the time interval considered) made and from the corresponding property of the maximum of two functions \cite{demmal} when this assumption is satisfied. 

Let $\overline \psi_1$ be a limit point of the sequence $\psi_{1_n}$. Then by the directional derivative definition we have 
$$\frac{\partial(-\ell(\overline \psi_1, x, z))}{\partial d} = -d_2 \overline \psi_1 + a_1 d_{1,1} \overline \psi_1 + a_2 d_{1,2} \overline \psi_1 + b \max\left\{ \left \langle \overline \psi_1 \frac{\partial f_1(x)} {\partial x}, d_1 \right \rangle, \left\langle \overline \psi_1 \frac{\partial f_2(x)} {\partial x}, d_1 \right \rangle \right\}.$$ From last two relations one obtains that condition iv) is fulfilled.

Thus, the function $h(x_1, x_2, z_1)$ satisfies conditions i)-iv), so it is differentiable in directions at the point $(x_1, x_2, z_1)$ \cite{BonnansShapiro}, and its derivative in the direction $d$ at this point is expressed by the formula
$$
\frac{\partial h(x_1, x_2, z_1)}{\partial d} = \sup_{\psi_1 \in \mathcal S (x_1, x_2, z_1)} \frac{\partial \ell(\psi_1, x_1, x_2, z_1)}{\partial d}
$$
where $\mathcal S (x_1, x_2, z_1) = \mathrm{arg} \max_{\psi_1 \in S_1} \ell(\psi_1, x_1, x_2, z_1)$. However, as shown above, in the problem considered the set $\mathcal S (x_1, x_2, z_1)$ consists of the only element $\psi_1^*(x_1, x_2, z_1)$, hence
$$
\frac{\partial h(x_1, x_2, z_1)}{\partial d} = \frac{\partial \ell(\psi_1^*(x_1, x_2, z_1), x_1, x_2, z_1)}{\partial d}.
$$

Finally, recall that by the directional derivative definition one has the equality
$$\frac{\partial (-\ell(\psi_1^*, x, z))}{\partial d} = -d_2 \psi_1^* + a_1 d_{1,1} \psi_1^* + a_2 d_{1,2} \psi_1^* + b \max\left\{ \left \langle \psi_1^* \frac{\partial f_1(x)} {\partial x}, d_1 \right \rangle, \left\langle \psi_1^* \frac{\partial f_2(x)} {\partial x}, d_1 \right \rangle \right\}$$
where we have put $\psi^*_1 := \psi_1^*(x_1, x_2, z_1)$.

From last two expressions we finally obtain that the function $h(x_1, x_2, z_1)$ is superdifferentiable at the point $(x_1, x_2, z_1)'$ but it is also positive in the case considered, so the function $h^2(x_1, x_2, z_1)$ is superdifferentiable at the point $(x_1, x_2, z_1)'$ as well as the square of a superdifferentiable positive function (see \cite{demrub}).

b) In the case $h_1 (x, z) = 0$ it is obvious that the function $h^2(x_1, x_2, z_1)$ is differentiable at the point $(x_1, x_2, z_1)'$ and its gradient vanishes at this point.

Now explore the differential properties of the functional $\varphi(x,z)$. For this, first give the following formulas for calculating the superdifferential $\overline \partial h^2(x, z)$ at the point $(x, z)$. At $i = \overline{1,n}$ one has
\begin{equation}
\label{5.5}
\overline \partial \left( {\textstyle \frac{1}{2}} \, h_i^2(x, z) \right) = h_i(x,z)\bigg(\psi^*_i {\bf e_{i+n}} - \psi^*_i [A'_i, {\bf 0_n}] + \sum_{j=1}^{r} \overline \partial \left(-\overline a_{ij} \max\big\{f_{i, j_1}(x) \psi_i^*, \dots, f_{i, j_{k(j)}}(x) \psi_i^*\big\} \right)\bigg)
\end{equation}
and at $j = \overline{1,r}$ we have
\begin{equation}
\label{5.7} \overline \partial \left(-\overline a_{ij} \max\big\{f_{i, j_1}(x) \psi_i^*, \dots, f_{i, j_{k(j)}}(x) \psi_i^*\big\} \right) = \overline a_{ij} \ \mathrm{co} \left \{ 
\left[ \psi_i^* \frac{\partial f_{i, j_p}(x)}{\partial x}, \bf{0_n} \right] 
\right \}, \quad j_p \in R_{ij}(x),
\end{equation}
$$
R_{ij}(x) = \Big\{j_p \in \{j_1, \dots, j_{k(j)}\} \ \big| \ f_{i, j_p}(x) \psi_i^* =  \max\big\{f_{i, j_1}(x) \psi_i^*, \dots, f_{i, j_{k(j)}}(x) \psi_i^*\big\} \Big\}.
$$

Calculate the derivative of the function ${\textstyle \frac{1}{2}} \, h_i^2(x, z) $ in the direction $g = [g_1, g_2] \in R^n \times R^n$. By virtue of formulas (\ref{5.5}), (\ref{5.7}) and superdifferential calculus rules \cite{demmal} we have
\begin{equation}
\label{100} 
 \frac {\partial\left({\textstyle \frac{1}{2}} \, h_i^2(x, z) \right) }{\partial g} = h_i(x,z)\bigg(\psi^*_i {{g_2}_i} - \left \langle \psi^*_i A_i,  g_1 \right \rangle + \sum_{j=1}^{r}  \min_{j_p \in R_{ij}(x)} \left\langle - \overline a_{ij} \psi_i^* \frac{\partial f_{i, j_p}(x)}{\partial x}, g_1 \right \rangle \bigg).
\end{equation}

Show that the functional $\varphi(x,z)$ is superdifferentiable and that its superdifferential is determined by the corresponding integrand superdifferential. 

\begin{thrm}
Let the interval $[0, T]$ may be divided into a finite number of intervals, in every of which one (several) of the functions $\big\{f_{i, j_1}(x) \psi_i, \dots, f_{i, j_{k(j)}}(x) \psi_i\big\}$, $i = \overline{1,n}$, $j = \overline{1,r}$, is (are) active. Then the functional
 $
 \varphi(x, z) 
 $
is superdifferentiable, i. e. 
\begin{equation}
\label{5.1} 
\frac{\partial \varphi(x, z)}{\partial g} = \lim_{\alpha \downarrow 0} \frac{1}{\alpha} \big(\varphi(x+\alpha g_1, z+\alpha g_2) - \varphi(x, z)\big) = \min_{w \in \overline \partial \varphi(x, z)} \int_0^T \langle w(t), g(t) \rangle dt \end{equation}
where $g = [g_1, g_2] \in C_n[0,T] \times P_n[0,T]$ and the set $\overline \partial \varphi(x, z)$ is defined as follows
\begin{equation}
\label{5.2'} 
\overline \partial \varphi(x, z) = \Big\{ w = [w_1, w_2] \in L_{n}^\infty[0, T] \times L_{n}^\infty[0, T] \ \big|
 \end{equation}
$$ [w_1(t), w_2(t)] \in \overline \partial \left( {\textstyle \frac{1}{2}} \, h^2(x(t), z(t)) \right) \quad for \ a. e. \ t \in [0, T] \Big\}.$$
\end{thrm}
\begin{proof}
In accordance with definition (\ref{0.6}) of a superdifferentiable functional, in order to prove the theorem, one has to check that

1) the derivative of the functional $\varphi(x, z)$ in the direction $g$ is actually of form (\ref{5.1}); 

2) herewith, the set $\overline \partial \varphi(x, z)$ is convex and weakly* compact subset of the space \linebreak $\left( C_n[0, T] \times P_n[0, T], || \cdot ||_{L^2_n [0, T] \times L^2_n [0, T]} \right)^*$. 

Let us prove statement 1). 

At first show that the following relations are true.

\begin{equation}
\label{6.16} 
\lim_{\alpha \downarrow 0}\frac{1}{\alpha} \bigg| \varphi(x+\alpha g_1, z+\alpha g_2) - \varphi(x,z) - \int_0^T  \min_{[w_1, w_2] \in \overline \partial \left(\frac{1}{2} h^2(x,z)\right)} \big( \langle w_1(t), \alpha g_1(t) \rangle + \langle w_2(t), \alpha g_2(t) \rangle \big) dt \bigg| = 0
 \end{equation} 
 or
 $$
\lim_{\alpha \downarrow 0}\frac{1}{\alpha} \bigg| \varphi(x+\alpha g_1, z+\alpha g_2) - \varphi(x,z) - \int_0^T   \frac {\partial\left({\textstyle \frac{1}{2}} \, h^2(x(t), z(t)) \right) }{\partial g}  dt \bigg| = 0.
$$
 
 
 Denote at $t \in [0, T]$
 \begin{equation}
\label{16'} 
\Phi(t, \alpha) = \frac{1}{\alpha} \Big ( h^2(x(t)+\alpha g_1(t), z(t)+\alpha g_2(t)) - h^2(x(t), z(t)) \Big) - 
 \end{equation} 
 $$ - \min_{[w_1, w_2] \in \overline \partial h^2(x,z)} \big( \langle w_1(t), g_1(t) \rangle + \langle w_2(t), g_2(t) \rangle \big) .$$

Our aim is to prove relation (\ref{6.16}) via Lebesgue's dominated convergence theorem applied to the function $\Phi(t, \alpha)$ (at $\alpha \downarrow 0$).

At first note that by superdifferential definition (\ref{0.1}) and by the superdifferentiability of the function $h^2(x,z)$ (proved at the beginning of this section) for each $t \in [0, T]$ we have $\Phi(t, \alpha) \rightarrow 0$ when $\alpha \downarrow 0$.  

In the following two paragraphs we show that for every $\alpha > 0$ one has $\Phi(t, \alpha) \in L_1^{\infty}[0,T]$.

Insofar as $x$, $g_1 \in C_n[0,T]$, $z$, $g_2 \in P_n[0,T]$ and the function $h^2(x,z)$ is continuous in its variables due to its structure \cite{demmal}, we obtain that for each $\alpha > 0$ the functions $t~\rightarrow~h^2(x(t),z(t))$ and $t \rightarrow h^2(x(t)+\alpha g_1(t)$, $z(t)+ \alpha g_2(t))$ belong to the space $L^\infty_1 [0,T]$.

Due to the upper semicontinuity of a superdifferential mapping \cite{demvas} and the structure of the superdifferential $\overline \partial h^2(x, z)$, it is easy to check that the mapping $t \rightarrow \overline \partial h^2(x(t), z(t))$ is upper semicontinuous, so it is measurable (see \cite{Blagodatskih}). Then due to the continuity of the function $g_1(t)$, the piecewise continuity of the function $g_2(t)$ and due to the continuity of the scalar product in its variables we obtain that the mapping 
\begin{equation}
\label{6.18} t \rightarrow \min_{[w_1, w_2] \in \overline \partial h^2(x(t),z(t))} \big( \langle w_1(t), g_1(t) \rangle + \langle w_2(t), g_2(t) \rangle \big) \end{equation} is upper semicontinuous \cite{aubenfr}, and then is also measurable \cite{Blagodatskih}. While proving statement 2) it will be shown that under the assumptions made the set $\overline \partial h^2(x, z)$ is bounded uniformly in $t \in [0, T]$, then by the continuity of the function $g_1(t)$ and the piecewise continuity of the function $g_2(t)$ it is easy to check that mapping (\ref{6.18}) is also bounded uniformly in $t \in [0, T]$. So we finally have that mapping (\ref{6.18}) belongs to the space $L^\infty_1 [0,T]$. 

Now we prove that the function {$\Phi(t, \alpha)$} is dominated by some integrable function on $[0, T]$ for all sufficiently small $\alpha > 0$. As is shown in the previous paragraph, the second summand in (\ref{16'}) is an integrable function. So it remains to consider for sufficiently small $\alpha > 0$ the first summand in (\ref{16'}). From the mean value theorem \cite{makela} (applied to the superdifferential) at each $t \in [0, T]$ and at each $\alpha > 0$ one has
$$ 
\frac{1}{\alpha} \Big ( h^2(x(t)+\alpha g_1(t), z(t)+\alpha g_2(t)) - h^2(x(t), z(t)) \Big) \in \overline \partial h^2\big(v_1(\alpha, t), v_2(\alpha_2, t) \big) g(t)
$$
where $v_1(\alpha, t) \in \mathrm{co}\big\{x(t), x(t) + \alpha g_1(t)\big\}$, $v_2(\alpha, t) \in \mathrm{co}\big\{z(t), z(t) + \alpha g_2(t)\big\}$. While proving statement 2) it will be shown that under the assumptions made the set $\overline \partial h^2(x, z)$ is bounded uniformly in $t \in [0, T]$, then by the continuity of the function $g_1(t)$ and the piecewise continuity of the function $g_2(t)$ the first summand in (\ref{16'}) is dominated by a piecewise continuous function for all sufficiently small $\alpha > 0$.

In \cite{Fom_opt_lett} it is shown that
\begin{equation}
\label{6.15}
\int_0^T \min_{[w_1, w_2] \in \overline \partial \left( \frac{1}{2} h^2(x,z)\right)}  \big( \langle w_1(t), g_1(t) \rangle + \langle w_2(t), g_2(t) \rangle \big) dt = \min_{w \in \overline \partial \varphi(x, z)} \int_0^T \langle w_1(t), g_1(t) \rangle + \langle w_2(t), g_2(t) \rangle dt. 
\end{equation}


From relations (\ref{6.16}), (\ref{6.15}) one obtains expression (\ref{5.1}). 

Prove statement 2): the corresponding proof may be found in \cite{Fominyh_sl_m}. 

The theorem is proved. 
\end{proof}

\begin{thrm} \label{th2'}
Let the interval $[0, T]$ may be divided into a finite number of intervals, in every of which one (several) of the functions $\big\{f_{i,j_1}(x) \psi_i, \dots, f_{i,j_{k(j)}}(x) \psi_i\big\}$, $i= \overline{1, n}$, $j = \overline{1,r}$, is (are) active. In order for the point $(x^{*}, z^*)$ to minimize the functional $I(x, z)$, it is necessary to have 
\begin{equation}
\label{5.444}
\overline \partial I (x^*(t), z^{*}(t)) = \{\bf 0_{2n}\}
\end{equation}
at each $t \in [0, T]$, where the expression for the superdifferential $\overline \partial I (x, z)$ is given by formula (\ref{4.3}) (where we take formula (\ref{5.2'}) for the superdifferential $\overline \partial \varphi(x,z)$). If one has $I(x^*, z^{*})$, then condition (\ref{5.444}) is also sufficient.
\end{thrm}

\section{Constructing the superdifferential descent direction \\ of the functional $I(x, z)$} \label{sc7}
In this section we consider only the points $(x, z)$ which do not satisfy minimum condition in Theorem \ref{th2'}. Our aim here is to find the superdifferential (or the steepest) descent direction of the functional $I(x, z)$ at the point $(x, z)$. Denote this direction $G(x, z)$. Herewith $G = [G_1, G_2] \in L^2_n[0, T] \times L^2_n[0, T]$. In order to construct the vector-function $G(x, z)$, consider the problem
\begin{equation}
\label{6.111}
\max_{w \in \overline{\partial} I(x, z)} ||w||^2_{L^2_{n}[0,T] \times L^2_{n}[0, T]}  = \max_{w \in \overline{\partial} I(x, z) } \int_0^T w^2(t)  dt. 
\end{equation}
Denote $\overline{w}$ the solution of this problem (below we will see that such a solution exists). The vector-function $\overline{w}$, of course, depends on the point $(x, z)$ but we omit this dependence in the notation for brevity. Then one can check that the vector-function 
\begin{equation}
\label{5.111} G\big(x(t), z(t), t\big) = -\frac{\overline{w}\big(x(t), z(t), t\big)}{||\overline{w}||_{L^2_{2n}[0,T] }}
\end{equation}
 is a superdifferential descent direction of the functional $I(x, z)$ at the point $(x, z)$ (cf. (\ref{6.66}), (\ref{6.6})). Recall that we are seeking the direction $G(x, z)$  in the case when the point $(x, z)$ does not satisfy minimum condition (\ref{5.444}), so $||\overline{w}||_{L_{2n}^{2}[0,T]} > 0$.    

Note that we have the equalities
$$\frac{\partial I(x, z)}{\partial G(x, z)} = \min_{w \in \partial \overline I(x, z)} \int_0^T \left \langle w(t), G(x(t), z(t), t) \right \rangle dt =\min_{w \in \overline \partial I(x, z)} \int_0^T \left \langle w(t), \frac{-\overline{w}(t)}{||\overline{w}||_{L_{2n}^{2}[0, T]}} \right \rangle dt =$$ 
$$= \frac{-1}{||\overline{w}||_{L_{2n}^{2}[0, T]}} \left( -\min_{w \in \overline \partial I(x, z)} \int_0^T \left \langle -w(t), {\overline{w}} \right \rangle dt \right)  = \frac{-1}{||\overline{w}||_{L_{2n}^{2}[0, T]}} \max_{w \in \partial \overline I(x, z)} \int_0^T \left \langle w(t), {\overline{w}(t)} \right \rangle dt = -||\overline{w}||_{L_{2n}^{2}[0,T]}$$
\newline
which, considering (\ref{0.6}) and the inequality $||\overline{w}||_{L_{2n}^{2}[0,T]} > 0$, implies 
\begin{equation}
\label{5.2}
I\left((x,z)+\alpha G(x,z)\right) < I(x,z) 
\end{equation}
for sufficiently small $\alpha > 0$.

It is easy to check that in this case the solution of problem (\ref{6.111}) is such a selector of the multivalued mapping $t \rightarrow \overline{\partial} I\big(x(t), z(t), t \big)$ that maximizes the distance from zero to the points of the set $\overline{\partial} I\big(x(t), z(t), t )$ at each time moment $t \in [0, T]$. In other words, to solve problem (\ref{6.111}) means to solve the following problem  
\begin{equation}
\label{5.3}
\max_{w(t) \in \overline{\partial} I(x(t), z(t), t ) } w^2(t)
\end{equation}
for each $t \in [0, T]$.
Actually, for every $t \in [0, T]$ we have the obvious inequality
$$
\max_{w \in \overline{\partial} I(x(t), z(t), t ) }  w^2(t) \geq w^2(t) 
$$
where $w(t)$ is a measurable selector of the mapping $t \rightarrow \overline{\partial} I\big(x(t), z(t), t \big)$ (from (\ref{4.2}) we have $w \in L^\infty_{2n} [0,T]$), then we obtain the inequality
$$
\int_0^T \max_{w \in \overline{\partial} I(x(t), z(t), t ) }  w^2(t) dt \geq \max_{w \in \overline{\partial} I(x, z) } \int_0^T w^2(t)  dt. 
$$
Insofar as for every $t \in [0, T]$ we have 
$$
\max_{w \in \overline{\partial} I(x(t), z(t), t ) }  w^2(t) \in \Big\{ w^2(t) \ \big| \ w(t) \in \overline{\partial} I(x(t), z(t), t ) \Big\}
$$
as the set $\overline{\partial} I\big(x(t), z(t), t \big) $ is closed and bounded at every fixed $t$ by definition of the superdifferential and the mapping $t \rightarrow \overline{\partial} I\big(x(t), z(t), t \big) $ is upper semicontinuous \cite {demvas} and besides, the norm (and its square) is continuous in its argument, then due to Filippov lemma \cite{Filippov2} there exists such a measurable selector $\overline{w}(t)$ of the mapping $t \rightarrow \overline{\partial} I\big(x(t), z(t), t \big) $ that for every $t \in [0, T]$ one obtains 
$$
\max_{w \in \overline{\partial} I(x(t), z(t), t ) }  w^2(t) = \overline{w}^2(t),
$$
so we have found the element $\overline{w}$ of the set $\overline{\partial} I(x, z) $ which brings the equality to the previous inequality. Hence, finally we obtain
$$
\int_0^T \max_{w \in \overline{\partial} I(x(t), z(t), t ) }  w^2(t) dt  = \max_{w \in \overline{\partial} I(x, z) } \int_0^T w^2(t)  dt. 
$$

Problem (\ref{5.3}) at each fixed $t \in [0, T]$ is a finite-dimensional problem of finding the maximal distance from zero to the points of the convex compact (the superdifferential). This problem can be effectively solved; the next paragraph describes its solution. In practice one makes a (uniform) partition of the interval $[0, T]$, and this problem is being solved for every point of the partition, i. e. one has to calculate $G\big(x(t_i), z(t_i), t_i \big)$ where $t_i \in [0, T]$, $i = \overline{1, N}$, are the points of discretization (see notation in Lemma \ref{lm1} below). Under some natural additional assumption Lemma \ref{lm1} below guarantees that the vector-function obtained with the help of piecewise linear interpolation of the superdifferential descent directions evaluated at every point of such a partition of the interval $[0, T]$ converges in the space $L^2_{2n}[0,T]$ to the vector-function $G\big(x(t), z(t), t\big)$ sought when the discretization rank tends to infinity.    

As noted in the previous paragraph, in order to construct the superdifferential descent direction of the functional $I(x, z)$ at the point $(x, z)$, it is required to find the maximal distance from zero to the points of the functional $I (x(t), z(t))$ superdifferential at each moment of time of a (uniform) partition of the interval~$[0, T]$. From formula (\ref{4.3}) (see also (\ref{5.5})) we see that the superdifferential $\overline \partial I (x(t), z(t))$ is a convex polyhedron $P(t) \subset~R^{2n}$. Herewith, of course, the set $P(t)$ depends on the point $(x, z)$. We will omit this dependence in the notation in this paragraph for simplicity. It is clear that in this case it is sufficient to go over all the vertexes $p_j(t)$, $j = \overline{1,s}$ (here $s$ is a number of vertexes of the polyhedron $P(t)$) and choose among the values $||p_j(t)||_{R^{2n}}$ the greatest one. Denote the corresponding vertex ${p}_{\overline j}(t)$ ($\overline j \in \{1, \dots, s\}$), and if there are several vertexes on which the maximal norm-value is achieved, then choose any of them. Finally, put $\overline w(t) = p_{\overline j}(t)$.

In work \cite{Fominyh_sl_m} one lemma is  given which, on the one hand, has rather natural for applications conditions and, on the other hand, guarantees that the function $L(t)$ obtained via piecewise linear interpolation of the function $v \in L_1^{\infty}[0, T]$ sought converges to this function in the space $L_1^2[0, T]$ when the rank of a (uniform) partition of the interval $[0, T]$ tends to infinity.   

\begin{lmm}\label{lm1}
 Let the function $v \in L_1^{\infty}[0, T]$ satisfy the following condition: for every $\overline{\delta} > 0$ the function $v(t)$ is piecewise continuous on the set $[0, T]$ with the exception of only the finite number of the intervals \linebreak $\big(\overline t_1(\overline{\delta}), \overline t_2(\overline{\delta})\big),$ $\dots$, $\big(\overline t_{r}(\overline{\delta}), \overline t_{r+1}(\overline{\delta})\big)$ whose union length does not exceed the value $\overline{\delta}$.  
 
 Choose the (uniform) finite splitting $t_1 = 0, t_2, \dots, t_{N-1}, t_N = T$ of the interval $[0, T]$ and calculate the values $v(t_i)$, $i = \overline{1, N}$, at these points. Let $L(t)$ be the function obtained via piecewise linear interpolation with the nodes $(t_i, v(t_i))$, $i = \overline{1, N}$. Then for every $\varepsilon > 0$ there exists such a number $\overline{N}(\varepsilon)$ that for every $N > \overline{N}(\varepsilon)$ one has $||L - v||^2_{L_1^2[0,T]} \leq \varepsilon$. 
 \end{lmm}
 
 At a qualitative level Lemma \ref{lm1} condition means that the sought function does not have ``too many'' points of discontinuity on the interval $[0, T]$. If this condition is satisfied for the vector-function $\overline w (t)$ (what it natural for applications), then this lemma justifies the approximation of the vector-function $\overline w(t)$ and hence, the approximation of the vector-function $\overline G\big( x(t), z(t), t \big)$, by the values $\overline w(t_i)$, $i = \overline{1,N}$, at the separate points of discretization implemented as described above.
 
\section{On a method for finding the stationary points of the functional $I(x, z)$} \label{sc8}
 Once the steepest (the superdifferential) descent direction has been constructed (see the previous section), one can apply some methods (based on using this direction) of nonsmooth optimization in order to find stationary points of the functional $I(x, z)$.

The simplest steepest (superdifferential) descent algorithm is used for numerical simulations of the paper.  In order to convey the main ideas, we firstly consider the convergence of this method for an analogous problem in a finite-dimensional case (which is of an independent interest); and then turn to a more general problem considered in this paper.


First consider the problem of minimization of a function which is a minimum of the finite number of continuously differentiable functions. So let 
$$
\varphi(x) = \min_{i = \overline{1, M}} f_i(x),
$$
where $f_i(x)$, $i = \overline{1, M}$, are continuously differentiable functions on $R^n$.  

It is known \cite{demmal} that the function $\varphi(x)$ is differentiable at every point $x_0 \in R^n$ in any direction $g \in R^n$, $||g||_{R^n} = 1$, and
\begin{equation}
\label{6.2}
\frac{\partial \varphi(x_0)}{\partial g} = \min_{i \in R(x_0)} \left \langle \frac{\partial f_i(x_0)}{\partial x}, g \right \rangle = \min_{w \in \overline\partial\varphi(x_0)} \langle w, g \rangle,
\end{equation}
$$R(x_0) = \{ i \in \{1,\dots,M\} \ | \ f_i(x_0) = \varphi(x_0) \},$$
$$\overline\partial\varphi(x_0) = \mathrm{co} \left\{ \frac{\partial f_i(x_0)}{\partial x}, \ i \in R(x_0) \right\}.$$

It is also known \cite{demmal} that in order for the point $x^* \in R^n$ to minimize the function $\varphi(x)$, it is necessary to have 
\begin{equation}
\label{6.3}
\min_{||g||_{R^n}=1} \min_{i \in R(x^*)} \left \langle \frac{\partial f_i(x^*)}{\partial x}, g \right \rangle \geq 0.
\end{equation}
Denote 
\begin{equation}
\label{6.4}
\Psi(x) = \min_{||g||_{R^n}=1}  \min_{i \in R(x)} \left \langle \frac{\partial f_i(x)}{\partial x}, g \right \rangle.
\end{equation}
Then necessary minimum condition (\ref{6.3}) of the function $\varphi(x)$ at the point $x^*$ may be rewritten as the inequality $\Psi(x^*) \geq 0$. 

From formulas (\ref{6.2}), (\ref{6.4}) we see that by definition the vector $g_0 \in R^n$, $||g_0||_{R^n} = 1$, is the steepest descent direction of the function $\varphi(x)$ at the point $x^0$ if we have
\begin{equation}
\label{6.5} \min_{i \in R(x_0)} \left \langle \frac{\partial f_i(x_0)}{\partial x}, g_0 \right \rangle = \Psi(x_0).
\end{equation}

In terms of superdifferential (see (\ref{6.2})) a steepest descent direction \cite{demvas} is the vector 
\begin{equation}
\label{6.66} g_0 = -\frac{w_0}{||w_0||_{R^n}} \end{equation} where $w_0$ is a solution of the problem
\begin{equation}
\label{6.6} \max_{w \in \overline\partial\varphi(x_0)} ||w|| = ||w_0||.
\end{equation}

Apply the superdifferential (the steepest) descent method to the function $\varphi(x)$ minimization. Describe the algorithm as applied to the function and the space under consideration. Fix an arbitrary initial point $x_1 \in R^n$. Suppose that the set $$lev_{\varphi(x_1)}\varphi(\cdot) = \{ x \in R^n \ | \ \varphi(x) \leq \varphi(x_1)\}$$ is bounded (due to the arbitrariness of the initial point, in fact, one must assume that the set $lev_{\varphi(x_1)}\varphi(\cdot)$ is bounded for every initial point taken). Due to the function $\varphi(x)$ continuity \cite{demvas} the set $lev_{\varphi(x_1)}\varphi(\cdot)$ is also closed. Let the point $x_k \in R^n$ be already constructed. If $\Psi(x_k) \geq 0$ (in practice we check that this condition is satisfied only with some fixed accuracy $\overline \varepsilon$, i. e. $\Psi(x_k) \geq -\overline \varepsilon$), then the point $x_k$ is a stationary point of the function $\varphi(x)$ and the process terminates. Otherwise, construct the next point according to the rule
$$x_{k+1} = x_k + \alpha_k g_k$$
where the vector $g_k$ is a steepest (superdifferential) descent of the function $\varphi(x)$ at
the point $x_k$ (see (\ref{6.66}), (\ref{6.6})) and the value $\alpha_k$ is a solution of the following one-dimensional minimization problem
$$\min_{\alpha \geq 0} \varphi(x_k + \alpha g_k) =  \varphi(x_k + \alpha_k g_k). $$
Then $\varphi(x_{k+1}) < \varphi(x_k)$. 

\begin{thrm}
Under the assumptions made one has the inequality
\begin{equation} 
\label{7_new_3}
\underline \lim_{k \rightarrow \infty} \Psi(x_k) \geq 0
\end{equation}
for the sequence built according to the algorithm above.
\end{thrm}
\begin{proof}
The following proof is partially based on the ideas of an analogous one in book \cite{demmal}.

Assume the contrary. Then there exist such a subsequence $\{x_{k_j}\}_{j=1}^{\infty}$ and such a number $b > 0$ that for each $j  \in \mathcal N$ we have the inequality
\begin{equation}
\label{7_new_1}
\Psi(x_{k_j}) \leq -b.
\end{equation} 
Note that for every $j \in \mathcal N$ the set $R(x_{k_j})$ is nonempty.

From the steepest descent direction definition (see (\ref{6.5})) it follows that for every $j \in \mathcal N$ there exists an index $\overline{i} \in R(x_{k_j})$ such that for each $\alpha > 0$ the relation
\begin{equation}
\label{6.7} f_{\overline i}(x_{k_j} + \alpha g_{k_j}) = f_{\overline i}(x_{k_j}) + \alpha \Psi(x_{k_j}) + o_{\overline i}(x_{k_j}, g_{k_j}, \alpha) \end{equation}
holds true.

Recall that $f_{\overline i}(x_{k_j}) = \varphi(x_{k_j})$ for each index $\overline i$ from the set $R(x_{k_j})$ by this index set definition. Then from inequality (\ref{7_new_1}) since $f_{\overline i}(x)$ is a continuously differentiable function by assumption, there exists \cite{demmal} such $\overline \alpha > 0$, which does not depend on the number $k_j$, such that for $\overline i \in R(x_{k_j})$ satisfying (\ref{6.7}) and for each $\alpha \in (0, \overline \alpha]$ one has the inequality
$$f_{\overline i}(x_{k_j} + \alpha g_{k_j}) \leq \varphi(x_{k_j}) - \frac{1}{2} \alpha b.$$

Using the fact that $R(x_{k_j}) \neq \emptyset$ for all $j \in \mathcal N$, we finally have
  \begin{equation}
\label{6.8} \varphi(x_{k_j} + \overline \alpha g_{k_j}) \leq \varphi(x_{k_j}) - \frac{1}{2} \overline \alpha b =  \varphi(x_{k_j}) - \frac{1}{2} \beta\end{equation}
  uniformly in $j \in \mathcal N$.
  
This inequality leads to a contradiction. Indeed, the sequence $\{\varphi(x_{k})\}_{k=1}^{\infty}$ is monotonically decreasing and bounded below by the number $\min\limits_{x \in lev_{\varphi(x_{(1)})} \varphi(\cdot)} \varphi(x)$, hence it has a limit:
  $$\{\varphi(x_{k})\} \rightarrow \varphi^* \ \mathrm{at} \ k \rightarrow \infty, $$
  herewith, at each $k \in \mathcal N$ one has
\begin{equation}
\label{6.9} 
\varphi(x_k) \geq \varphi^*.
\end{equation}

Now choose such a big number $\overline j$ that 
$$\varphi(x_{k_{\overline j}}) < \varphi^* + \frac{1}{2} \beta.$$
Due to (\ref{6.8}) we have
$$\varphi(x_{k_{\overline j} + 1}) \leq \varphi(x_{k_{\overline j}} + \overline \alpha g_{k_{\overline j}}) \leq \varphi^* - \frac{1}{2} \beta$$
what contradicts (\ref{6.9}).
\end{proof}

Now turn back to the problem of functional $I(x,z)$ minimization. Denote (see formulas (\ref{6.111}), (\ref{5.111}))
\begin{equation}
\label{6.10}
\Psi(x, z) = \min_{||g||_{L_n^2[0, T] \times L_n^2[0,T]}=1}  \frac{\partial I(x, z)}{\partial g} = \frac{\partial I(x, z)}{\partial G} .
\end{equation}
Then necessary minimum condition of the functional $I(x, z)$ at the point $(x^*, z^*)$ may be written \cite{demvar} as the inequality $\Psi(x^*, z^*) \geq 0$. 

Apply the superdifferential (the steepest) descent method to the functional $I(x, z)$ minimization. Describe the algorithm as applied to the functional and the space under consideration. Fix an arbitrary initial point $(x_{(1)}, z_{(1)}) \in C_n[0, T] \times P_n[0, T]$. Suppose that the set $$lev_{I(x_{(1)}, z_{(1)})}I(\cdot, \cdot) = \left\{ (x, z) \in C_n[0,T] \times P_n[0,T] \ | \ I(x, z) \leq I(x_{(1)}, z_{(1)})\right\}$$ is bounded in $L_n^2[0, T] \times L_n^2[0,T]$-norm (due to the arbitrariness of the initial point, in fact, one must assume that the set $lev_{I(x_{(1)}, z_{(1)})}I(\cdot, \cdot)$ is bounded for every initial point taken). Let the point $(x_{(k)}, z_{(k)}) \in C_n[0, T] \times P_n[0, T]$ be already constructed. If $\Psi(x_{(k)}, z_{(k)}) \geq 0$ (in practice we check that this condition is satisfied only with some fixed accuracy $\overline \varepsilon$, i. e. $\Psi(x_{(k)}, z_{(k)}) \geq -\overline \varepsilon$) (in other words, if minimum condition  (\ref{5.444}) is satisfied (in practice with some fixed accuracy $\overline{\varepsilon}$, i. e. $||\overline{w}\big(x_k(t), z_k(t), t\big)||_{L_{2n}^2[0, T]} \leq \overline{\varepsilon}$)), then the point $(x_{(k)}, z_{(k)})$ is a stationary point of the functional $I(x, z)$ and the process terminates. Otherwise, construct the next point according to the following rule
$$
(x_{(k+1)}, z_{(k+1)}) = (x_{(k)}, z_{(k)}) + \alpha_{(k)} G\big(x_{(k)}, z_{(k)}\big)
$$
where the vector-function $G\big(x_{(k)}, z_{(k)}\big)$ is a superdifferential descent direction of the functional $I(x, z)$ at the point $(x_{(k)}, z_{(k)})$ (see (\ref{6.111}), (\ref{5.111})), and the value $\alpha_{(k)}$ is a solution of the following one-dimensional minimization problem
$$
\min_{\alpha \geq 0} I \Big( (x_{(k)}, z_{(k)}) + \alpha G\big(x_{(k)}, z_{(k)}\big) \Big) = I \Big( (x_{(k)}, z_{(k)}) + \alpha_{(k)} G\big(x_{(k)}, z_{(k)}\big) \Big). 
$$
Then according to (\ref{5.2}) one has $$
I \big(x_{(k+1)}, z_{(k+1)}\big) < I \big(x_{(k)}, z_{(k)}\big). $$

Introduce now the set family $\mathcal{I}$. At first define the functional $I_q$, $q = \overline{1, \left(\prod_{j=1}^r {k(j)}\right)^n}$ as follows. Its integrand is the same, as the functional $I$ one, but the maximum function $\max\big\{f_{i, j_1}(x) \psi_i, \dots, f_{i, j_{k(j)}}(x) \psi_i\big\}$, $j = \overline{1,r}$, is substituted for each $i = \overline{1,n}$ by only the one of the functions $f_{i, j_1} \psi_i, \dots, f_{i, j_{k(j)}} \psi_i$, $j \in \{1, \dots, r\}$. Let the family $\mathcal{I}$ consist of the sums of the integrals over the intervals of the time interval $[0, T]$ splitting for all possible such finite splittings. Herewith, the integrand of each summand in the sum taken is the same, as some functional $I_q$ one, $q \in \left\{1, \dots, \left(\prod_{j=1}^r {k(j)}\right)^n \right\}$.

Let for every point constructed by the method described the following assumption is valid: the interval $[0, T]$ may be divided into a finite number of intervals, in every of which for each $i =\overline{1, n}$ either $h_i(x_{(k)},z_{(k)}) =0$, or one (several) of the functions $\left\langle - \overline a_{ij} \psi_i^* \frac{\partial f_{i, j_p}(x_{(k)})}{\partial x}, G_1(x_{(k)},z_{(k)}) \right \rangle$, $j = \overline{1,r}$, $p = \overline{1,k(j)}$, is (are) active.

Let us illustrate this assumption by an example. Consider the following simplest functional (whose structure, however, preserves the basic features of the general case)
$$\int_0^1 \frac{1}{2} \mathrm{min}^2\{x(t)+1, -x(t)+1\} dt.$$

Consider the point $x_{(1)} = 0$ (i. e. $x_{(1)}(t) = 0$ for all $t \in [0, 1]$). In order to find the steepest (the superdifferential) descent direction of the functional at this point one has, according to the theory described, to minimize the directional derivative (calculated at this point), i. e. to find such a function $G \in L_1^2[0, T]$, $\displaystyle{\int_0^1 G^2(t)} dt = 1$, that minimizes the functional 
$$\int_0^1 \min\{x_{(1)}(t)+1, -x_{(1)}(t)+1\} \min\{g(t), -g(t)\} dt.$$
Here $g \in L_1^2[0, T]$, $\displaystyle{\int_0^1 g^2(t)} dt = 1$. Take 
$$G(t) =  \left\{
\begin{array}{ll}
-1, \ &\text{if} \ t \in [0, 0.5],  \\
1, \ &\text{if} \ t \in (0.5, 1]
\end{array}
\right.
$$
as one of obvious solutions. Herewith, $\Psi(x_{(1)}) = \displaystyle \int_0^1 -1 dt = -1$. 
We see that the assumption made is satisfied. Take 
$$\hat I(x) = \int_0^{0.5} \frac{1}{2} (x(t)+1)^2 dt + \int_{0.5}^1 \frac{1}{2} (-x(t)+1)^2 dt.$$
Then we have
$$ \hat I\left( x_{(1)} + \alpha G(x_{(1)}) \right) = \hat I (x_{(1)})  + \alpha \Psi(x_{(1)}) + o(x_{(1)}, G(x_{(1)}), \alpha) $$
since $$\nabla \hat I(x_{(1)}) = \left\{
\begin{array}{ll}
1, \ &\text{if} \ t \in [0, 0.5],  \\
-1, \ &\text{if} \ t \in (0.5, 1],
\end{array}
\right.$$ 
i. e. $\nabla \hat I(x_{(1)}) G(x_{(1)}) = \displaystyle \int_0^1 -1 dt = -1$. It is obvious that $\hat I \in \mathcal{I}$. Note that with $\alpha_{(1)} = 1$ one gets the point 
$$x_{(2)}(t) =  \left\{
\begin{array}{ll}
-1, \ &\text{if} \ t \in [0, 0.5],  \\
1, \ &\text{if} \ t \in (0.5, 1],
\end{array}
\right.
$$
which delivers the global minimum to the functional considered.

For functionals from the family $\mathcal I$ we make the following additional assumption. Let there exists such a finite number $L$ that for every $\hat I \in \mathcal I$ and for all $\overline x, \overline z, \overline{\overline x}, \overline{\overline z}$ from a ball with the center in the origin and with some finite radius $R' + \hat \alpha$ (here $R' > \sup\limits_{(x,z) \in lev_{I(x_{(1)}, z_{(1)})}I(\cdot, \cdot)}||(x, z)||_{L^2_n[0,T] \times L^2_n[0,T] }$ and $\hat \alpha$ is some positive number) one has 
\begin{equation}
\label{6.100} ||\nabla \hat I(\overline x, \overline z) -  \nabla \hat I(\overline{\overline x}, \overline{\overline z})||_{L^2_n[0,T] \times L^2_n[0,T] } \leq L ||({\overline x}, {\overline z})' - (\overline{\overline x}, \overline{\overline z})' ||_{L^2_n[0,T] \times L^2_n[0,T]}. 
\end{equation}

\begin{rmrk} \label{rmrk3}
At first glance it may seem that the Lipschitz constant $L$ existence for all $\hat I \in \mathcal{I}$ simultaneously in the assumption made is too burdensome. However, if one remembers that on each of the finite number of the interval $[0, T]$ segments the functional $\hat I \in \mathcal I$ integrand coincides with the integrand of the functional $I_q$, $q \in \left\{1, \dots, \left(\prod_{j=1}^r {k(j)}\right)^n \right\}$, by construction (see the set $\mathcal I$ definition), then this assumption is natural if we suppose the Lipschitz-continuity of every of the gradients $\nabla I_q$, $q = \overline{1, \left(\prod_{j=1}^r {k(j)}\right)^n}$; and this gradient Lipschitz-continuity condition is a common assumption for justifying classical optimization methods for differentiable functionals.
\end{rmrk}

\begin{lmm}\label{lm2} Let condition (\ref{6.100}) be satisfied. Then for each functional $I \in \mathcal I$ and for all $(x,z) \in lev_{I(x_{(1)}, z_{(1)})}I(\cdot, \cdot)$, $G \in C_n[0, T] \times P_n[0, T]$, $||G||_ {L_n^2[0,T] \times L_n^2[0,T] } = 1$, $\alpha \in R$, $0 < \alpha \leq \hat \alpha$ one has the inequality
$$
I\left( (x, z) + \alpha G \right) \leq I(x,z) + \alpha \left \langle \nabla I(x,z), G \right \rangle + \alpha^2 \frac{L}{2}. 
$$
\end{lmm}

\begin{proof} The proof can be carried out with obvious modifications in a similar way as for the analogous statement in \cite{KantorovichAkilov}. \end{proof}

We suppose that during the method realization for each $ k \in \mathcal N$ one has $\alpha_{(k)} < \hat \alpha$ where $\hat \alpha$ is a number from Lemma \ref{lm2} (see also the assumption before Remark \ref{rmrk3}).

\begin{thrm}
Under the assumptions made one has the inequality
\begin{equation} 
\label{7_new_4}
\underline \lim_{k \rightarrow \infty} \Psi(x_k, z_k) \geq 0
\end{equation}
for the sequence built according to the rule above.
\end{thrm}
\begin{proof}
Assume the contrary. Then there exist such a subsequence $\{(x_{k_j},z_{k_j})\}_{j=1}^{\infty}$ and such a number $b > 0$ that for each $j  \in \mathcal N$ we have the inequality
\begin{equation}
\label{7_new_2}
\Psi(x_{k_j},z_{k_j}) \leq -b.
\end{equation}

From the steepest descent direction (see (\ref{6.10})) and the set $\mathcal I$ definitions (see also (\ref{100})) it follows that for every $j \in \mathcal N$ there exists a functional $\hat I$ from the family $\mathcal I$ such that for each $\alpha > 0$ the relation
\begin{equation}
\label{6.11} \hat I\left( (x_{k_j}, z_{k_j}) + \alpha G(x_{k_j}, z_{k_j}) \right) = \hat I (x_{k_j}, z_{k_j})  + \alpha \Psi(x_{k_j},z_{k_j}) + o\left( (x_{k_j}, z_{k_j}), G(x_{k_j}, z_{k_j}), \alpha \right) \end{equation}
holds; herewith,
$$\Psi(x_{k_j}, z_{k_j}) = \left \langle \nabla \hat I(x_{k_j}, z_{k_j}), G  (x_{k_j}, z_{k_j}) \right \rangle.$$

Recall that $\hat I (x_{k_j}, z_{k_j})  = I (x_{k_j}, z_{k_j}) $ for each functional $\hat I$ from the family $\mathcal I$ by this set definition. Then from inequality (\ref{7_new_2}) by virtue of Lemma \ref{lm2} there exists $\overline \alpha > 0$, which does not depend on the number $k_j$, such that for $\hat I \in \mathcal{I}$ satisfying (\ref{6.11}) and for each $\alpha \in (0, \overline \alpha]$ one has the inequality  
$$\hat I\left( (x_{k_j}, z_{k_j}) + \alpha G(x_{k_j}, z_{k_j}) \right)  \leq I(x_{k_j}, z_{k_j}) - \frac{1}{2} \alpha b .$$

Using the set $\mathcal I$ definition, we finally have
  \begin{equation}
\label{6.88} I\left( (x_{k_j}, z_{k_j}) + \alpha G(x_{k_j}, z_{k_j}) \right) \leq I(x_{k_j}, z_{k_j}) - \frac{1}{2} \overline \alpha b =  I(x_{k_j}, z_{k_j}) - \frac{1}{2} \beta \end{equation}
  uniformly in $j \in \mathcal N$.

This inequality leads to a contradiction. Indeed, the sequence $\left\{I(x_{(k)},z_{(k)})\right\}_{k=1}^{\infty}$ is monotonically decreasing and bounded below by zero (recall that the functional $I(x,z)$ is nonnegative by construction), hence, it has a limit:
  \begin{equation}
\label{7.100} \left\{I(x_{(k)},z_{(k)})\right\} \rightarrow I^* \ \mathrm{at} \ k \rightarrow \infty, \end{equation}
  herewith, at each $k \in \mathcal N$ one has $\{I(x_{(k)},z_{(k)})\} \geq I^*$.

Now choose such a big number $\overline j$ that 
$$I(x_{(k_{\overline j})},z_{(k_{\overline j})}) < I^* + \frac{1}{2} \beta.$$
Due to (\ref{6.88}) we have 
$$I(x_{(k_{\overline j}+1)},z_{(k_{\overline j}+1)}) \leq I(x_{(k_{\overline j})},z_{(k_{\overline j})}) + \overline \alpha G_{k_{\overline j}}) \leq I^* - \frac{1}{2} \beta$$
what contradicts (\ref{7.100}).
\end{proof}

\begin{rmrk}
It is easy to show that, in fact, in formulas (\ref{7_new_3}), (\ref{7_new_4}) the lower limit can be substituted by the ``ordinary'' limit and the inequality can be substituted by the equality in the cases considered.
\end{rmrk}

\section{Numerical examples} \label{sc9}
In this section we give two illustrative examples of the described algorithm implementation.  

\begin{xmpl} Consider the differential inclusion
$$\dot x_1 \in x_2 + [-2, 2] (|x_1| + |x_2|), \quad \dot x_2 \in x_1 + [-2, 2] (|x_1| + |x_2|)$$
on the time interval $[0, 1]$ with the boundary conditions
$$x_1(0) = 0, \, x_2(0) = 1, \quad x_1(1) = 0, \, x_2(1) = 2.$$
As the phase constraint surface put $$s(x) = x_1 = 0.$$ 

Take $(x_{(1)}, z_{(1)}) = (0, 0, 0, 0)'$ as the first approximation (i. e. all the functions considered are identically equal to zero on the interval $[0,1]$), then $I(x_{(1)}, z_{(1)}) = 1$. At the end of the process the discretization step was equal to $2 \times 10^{-1}$. Figure 1 illustrates the trajectories obtained. From the picture we see that the differential inclusion is satisfied. The trajectory lies on the surface required as well. The boundary values error doesn't exceed the magnitude $5 \times 10^{-3}$. In order to obtain such an accuracy 5 iterations has been required. The functional value on the trajectory obtained is approximately $3 \times 10^{-5}$.

\begin{figure*}[h!]
\begin{minipage}[h]{0.3\linewidth}
\center{\includegraphics[width=1\linewidth]{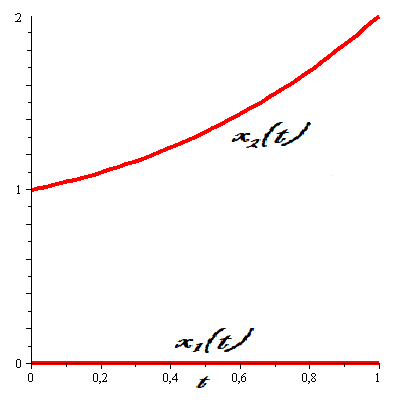} }
\end{minipage}
\hfill
\begin{minipage}[h]{0.3\linewidth}
\center{\includegraphics[width=1\linewidth]{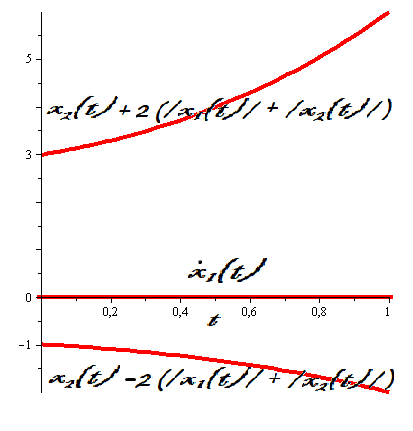} }
\end{minipage}
\hfill
\begin{minipage}[h]{0.3\linewidth}
\center{\includegraphics[width=1\linewidth]{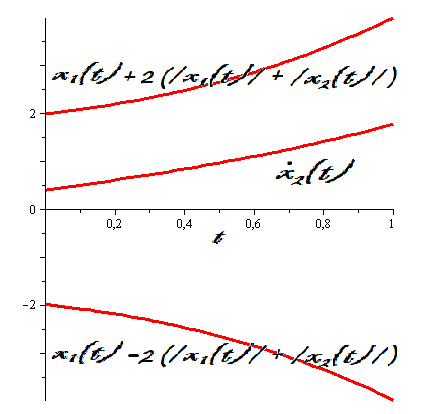} }
\end{minipage}
\caption{Example 1}
\end{figure*}
\end{xmpl}
 
\begin{xmpl} Consider the differential inclusion
$$\dot x_1 \in {\mathrm{co} \{x_1, x_2, x_3\}}, \quad \dot x_2 \in {\mathrm{co} \{-2 x_1, x_3\}}, \quad \dot x_3 = x_1+x_2$$
on the time interval $[0, 1]$ with the boundary conditions
$$x_1(0) = -1, \, x_2(0) = -1, \, x_3(0) = 2.5, \quad x_1(1) = 1, \, x_2(1) = 1, \, x_3(1) = 2.5.$$
As the phase constraint surface put $$s(x) = x_1-x_2 = 0.$$ 

Take $(x_{(1)}, z_{(1)}) = (0, 0, 0, 0, 0, 0)'$ as the first approximation (i. e. all the functions considered are identically equal to zero on the interval $[0,1]$), then $I(x_{(1)}, z_{(1)}) = 8.125$. At the end of the process the discretization step was equal to $10^{-1}$. Figure 2 illustrates the trajectories obtained. From the picture we see that the differential inclusion is satisfied. The trajectory lies on the surface required as well. The boundary values error doesn't exceed the magnitude $2 \times 10^{-3}$. In order to obtain such an accuracy 66 iterations has been required. The functional value on the trajectory obtained is approximately $2 \times 10^{-5}$.

Note that we have intentionally taken such a constraint surface in order to make the trajectories $x_1(t)$ and $x_2(t)$ coincide. So if on some iterations the trajectory already lies on the surface $s(x) = 0$ and $h_1(x, z) > 0$, then both of these trajectories would be active on some nonzero time interval (if $x_1(t) = x_2(t) > x_3(t)$ on this interval). That would lead to calculating the ``full-fledged'' superdifferential (not the one that degenerates into a single point) on these iterations.   

\begin{figure*}[h!]
\begin{minipage}[h]{0.3\linewidth}
\center{\includegraphics[width=1\linewidth]{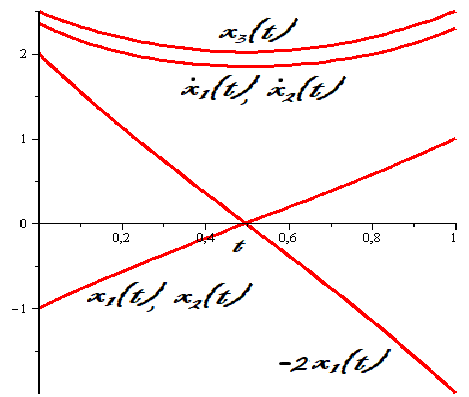} }
\end{minipage}
\hfill
\begin{minipage}[h]{0.3\linewidth}
\center{\includegraphics[width=1\linewidth]{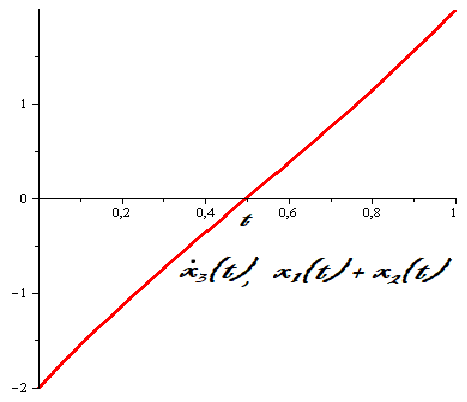} }
\end{minipage}
\caption{Example 2}
\end{figure*}
\end{xmpl}

{The main results of this paper (Sections \ref{sc6}, \ref{sc7}, \ref{sc8}, \ref{sc9}) were obtained in IPME RAS and supported by Russian Science Foundation (grant 20-71-10032).}

\end{document}